\documentclass{article}
\usepackage{amsfonts,amsmath}
\mathchardef\mhyphen="2D 
\usepackage[a4paper]{geometry}
\usepackage{microtype}
\usepackage{relsize}

\def\zo/{$0\mkern2mu\mhyphen1$}

\def\nn/{$n \times n$}

\title{A Near Proof of Weak Graph Positivity, A New Property of Regular Random Graphs}
\date{\today} 
\author{Paul Federbush\\
Department of Mathematics\\
University of Michigan\\
Ann Arbor, MI, 48109-1043}

\usepackage{amsthm} 
\newtheorem{thm}{Theorem}[section]
\newtheorem*{conj*}{Conjecture} 
\numberwithin{equation}{section} 
\DeclareMathOperator{\Prob}{Prob}
\DeclareMathOperator{\EE}{E}
\DeclareMathOperator{\OO}{\mathcal{O}}
\DeclareMathOperator{\PP}{P}
\newcommand{\sslash}{\mathbin{/\mkern-6mu/}} 

\begin{document}

\maketitle
\begin{abstract}
  One deals with $r$-regular bipartite graphs with $2n$ vertices.
  In a previous paper Butera, Pernici, and the author have introduced a quantity
  $d(i)$, a function of the number of $i$-matchings, and conjectured that as $n$
  goes to infinity the fraction of graphs that satisfy $\Delta^k d(i) \ge 0$,
  for all $k$ and $i$, approaches $1$.
  Here $\Delta$ is the finite difference operator. This conjecture we called the 'graph
  positivity conjecture'.
  In this paper it is formally shown that for each $i$ and $k$  the
  probability that $\Delta^k d(i) \ge 0$ goes to $1$ with $n$ going to infinity.
  We call this weaker result the 'weak graph positivity conjecture ( theorem )'.
  A formalism of Wanless as systematized by Pernici is central to this effort.
  Our result falls short of being a rigorous proof since we make a sweeping
  conjecture ( computer tested ), of which we so far have only a portion of the proof.

\end{abstract}

\section{Introduction}
We deal with $r$-regular bipartite graphs with $v = 2n$ vertices.
We let $m_i$ be the number of $i$-matchings.
In \cite{BFP}, Butera, Pernici, and I introduced the quantity $d(i)$, in eq.\
$(10)$ therein,
\begin{equation}\label{eq:1.1}
  d(i) \equiv \ln\biggl( \frac{m_i}{r^i} \biggr)
  - \ln\biggl( \frac{\overline{m}_i}{(v-1)^i} \biggr)
\end{equation}
where $\overline{m}_i$ is the number of $i$-matchings for the complete (not
bipartite complete) graph on the same vertices,
\begin{equation}\label{eq:1.2}
  \overline{m}_i = \frac{v!}{(v-2i)!\,i!\,2^i}
\end{equation}
We here have changed some of the notation from \cite{BFP} to agree with notation
in \cite{Per}.
We then considered $\Delta^k d(i)$ where $\Delta$ is the finite difference
operator, so
\begin{equation}\label{eq:1.3}
  \Delta d(i) = d(i+1) - d(i)
\end{equation}
A graph was defined to satisfy \textit{graph positivity} if all the meaningful
$\Delta^k d(i)$ were non-negative.
That is
\begin{equation}\label{eq:1.4}
  \Delta^k d(i) \ge 0
\end{equation}
for $k = 0,\ldots,v$ and $i = 0,\ldots,v-k$.
We made the conjecture, the 'graph positivity conjecture', supported by some computer evidence,
\begin{conj*}
  As $n$ goes to infinity the fraction of graphs that satisfy graph positivity
  approaches one.
\end{conj*}
In this paper we "formally" prove a weaker result, the 'weak graph positivity conjecture',
 in the same direction.
\begin{thm}\label{thm:1.1}
  For each $i$ and $k$, one has
  \begin{equation*}
    \Prob\bigl(\Delta^k d(i) \ge 0\bigr) \xrightarrow[n \to \infty]{} 1
  \end{equation*}
\end{thm}
The paper relies heavily on the work of Wanless, \cite{Wan}, and Pernici,
\cite{Per}, that gives a nice representation of the $m_i$. 
The restriction to bipartite graphs is mainly because this restriction is made in \cite{BFP}. In this paper the
bipartite nature appears used in two places. First, the number of vertices is assumed to be
even, and second, in eq.(3.5) the lower limit 4 is replaced by 3 if one does not assume
the graph is bipartite. 

As said in the abstract our result falls short of being a rigorous proof, since we make a sweeping
  conjecture ( computer tested )  presented in Section 10. We are assembling parts of the proof
  and hope to have a complete proof in the not too distant future.
  
  One should read the \textit{Valuable Observation} at the end of Section 4 to see why the question
  of convergence ( of series ) in this paper is a piece of cake.
  
  We hope the reader sees the beauty of some of the arguments in this paper, simple ideas building
  on one another.

\section{Idea of the proof}
Suppose we want $\Prob(x > y)$ to be large.
We have
\begin{equation}\label{eq:2.1}
  \Prob(x < y) = \Prob(e^x < e^y)
\end{equation}
Set
\begin{equation}\label{eq:2.2}
 (e^x - e^y) \equiv \alpha_0
\end{equation}
and
\begin{equation}\label{eq:2.3}
  \EE(e^x - e^y) \equiv\alpha
\end{equation}
We will want $ \alpha$ to be positive, and in the present work proving this
positivity will be a major component.
Let
\begin{equation}\label{eq:2.3a}
  \EE\bigl( (e^x - e^y)^2\bigr) - \alpha^2 \equiv \beta
\end{equation}
Then, assuming $\alpha > 0$,
\begin{equation}\label{eq:2.4}
  \EE\bigl( (e^x - e^y - \alpha)^2\bigr) \ge
   \alpha^2 \, \PP(e^x - e^y < 0)
\end{equation}
And so
\begin{equation}\label{eq:2.5}
  \Prob(e^x < e^y) \le \frac{\beta}{\alpha^2}
\end{equation}
In our problem $\beta$ and $\alpha$ will be functions of $n$, and we'll want
probability to go to zero with $n$ as $n$ goes to infinity.
\par We turn to the object of study, and perform some simple manipulations,
working from eq.\ \eqref{eq:1.1}
\begin{align}
  &\Prob\bigl( \Delta^k d(i) > 0 \bigr) = \Prob\Biggl( (-1)^k
  \sum_{\ell=0}^k (-1)^\ell \binom{k}{\ell} \, \ln \biggl(
  \frac{m_{i+\ell}}{\overline{m}_{i+\ell}} \cdot \frac{(v-1)^{i+\ell}}{r^{i+\ell}}
  \biggr) > 0 \Biggr)\label{eq:2.6}\\
  &\qquad= \Prob\Biggl( \sum_{\ell \in \mathcal{L}^+} \binom{k}{\ell} \,
  \ln \biggl( \frac{m_{i+\ell}}{\overline{m}_{i+\ell}} \cdot
  \frac{(v-1)^{i+\ell}}{r^{i+\ell}} \biggr) > \sum_{\ell \in \mathcal{L}^-}
  \binom{k}{\ell}\,\ln \biggl( \frac{m_{i+\ell}}{\overline{m}_{i+\ell}} \cdot
  \frac{(v-1)^{i+\ell}}{r^{i+\ell}} \biggr) \Biggr)\label{eq:2.7}
\end{align}
where $\mathcal{L}^+$ is the set of odd $\ell$, $0 \le \ell \le k$, if $k$ is
odd and is the set of even $\ell$, $0 \le \ell \le k$, if $k$ is even, and
$\mathcal{L}^-$ is defined vice versa.
\par Returning to the language of \eqref{eq:2.1}--\eqref{eq:2.5}, we set
\begin{align}
  x &\equiv \sum_{\ell \in \mathcal{L}^+} \binom{k}{\ell}\,\ln \biggl(
  \frac{m_{i+\ell}}{\overline{m}_{i+\ell}} \cdot
  \frac{(v-1)^{i+\ell}}{r^{i+\ell}} \biggr)\label{eq:2.8}\\
  y &\equiv \sum_{\ell \in \mathcal{L}^-} \binom{k}{\ell}\,\ln \biggl(
  \frac{m_{i+\ell}}{\overline{m}_{i+\ell}} \cdot
  \frac{(v-1)^{i+\ell}}{r^{i+\ell}} \biggr)\label{eq:2.9}
  \intertext{and so}
  e^x &= \prod_{\ell \in \mathcal{L}^+} \biggl(
  \frac{m_{i+\ell}}{\overline{m}_{i+\ell}} \cdot
  \frac{(v-1)^{i+\ell}}{r^{i+\ell}} \biggr)^{\binom{k}{\ell}}\label{eq:2.10}\\
  e^y &= \prod_{\ell \in \mathcal{L}^-} \biggl(
  \frac{m_{i+\ell}}{\overline{m}_{i+\ell}} \cdot
  \frac{(v-1)^{i+\ell}}{r^{i+\ell}} \biggr)^{\binom{k}{\ell}}\label{eq:2.11}
\end{align}
Throughout, the caveat "for $n$ large enough" is understood.
\section{The Work of Wanless and Pernici}\label{sec:3}
In \cite{Wan} Wanless developed a formalism to compute the $m_i$ of any regular
graph.
We here only give a flavor of this formalism, but present some of the consequences we
will use in this paper.
For each $i$ there are defined a set of graphs
$g_{i1},g_{i2},\ldots,g_{i n(i)}$.
Given a regular graph $g$, one computes for each $j$ the number of subgraphs of
$g$ isomorphic to $g_{ij}$, call this $g \sslash g_{ij}$.
Then $m_i$ for $g$ is determined by the $n(i)$ values of $g \sslash g_{ij}$.
We define $M_i$ to be the value of $m_i$ assigned to any graph with all $n(i)$
values of $g \sslash g_{ij}$ zero. Such graphs will exist only for large enough n.
Initially $M_i$ is defined only for such n. But it may be extended as a finite
polynomial in $\frac{1}{n}$ to all non-zero $n$.
$M_i = M_i(r,n)$ is an important object of study to us.
\par In \cite{Per} Pernici systematized the results of Wanless; we take a number
of formulae from this paper.
From eq.\ $(12)$ and eq.\ $(13)$ of \cite{Per} we write
\begin{align}
  M_j &= \frac{n^jr^j}{j!} (1+H_j) \label{eq:3.2}\\
  H_j &= \sum_{h=1}^{j-1} \frac{a_h(r,j)}{n^h} \label{eq:3.3}
\end{align}
$a_h$ is a polynomial of degree at most $2h$ in $j$.  We view $M_j$ and $H_j$ as formal
polynomials in $j$ amd $\frac{1}{n}$.
We will sometimes need values of $a_h$, as given by eq.\ $(18)$ and eq.\ $(45)$ of \cite{Per}.

We will use eq.\ $(16)$ and $(17)$ of \cite{Per}
\begin{align}
   [j^k n^{-h}]\,\ln \biggl( 1 + H_j \biggr) =
    [j^k n^{-h}]\,\ln \biggl( 1 + \sum_{s=1}^{j-1} \frac{a_s(r,j)}{n^s} \biggr) &=
  0, \qquad k \ge h+2 \label{eq:3.4}\\
  [j^{h+1} n^{-h}]\,\ln \biggl( 1 + H_j \biggr)=  [j^{h+1} n^{-h}]\,\ln \biggl( 1 + \sum_{s=1}^{j-1} \frac{a_s(r,j)}{n^s} \biggr)
  &= \frac{1}{(h+1)h} \biggl( \frac{1}{r^h} - 2 \biggr) \label{eq:3.5}
\end{align}
where $[j^k n^{-h}]$ in front of an expression picks out the coefficient of the
$\frac{j^k}{n^h}$ term, $c(k,h)$, in an expression $\sum_{\alpha,\beta}
c(\alpha,\beta) \frac{j^\alpha}{n^\beta}$. We are here working with the formal
polynomials with $r$ fixed.  Eq. (3.3) and (3.4) are special cases 
of the conjecture of Section 10. We already have a rigorous proof of (3.3) due to
Robin Chapman, \cite{5}. We have some good ideas toward the full proof.
\par We set $M_0 = 1$ and $M_s = 0$ if $s < 0$.
Then $m_j$ is recovered from $M_j$ by the formula
\begin{equation}\label{eq:3.6}
  m_j = \exp \biggl( \sum_{s \ge 4} \frac{\varepsilon_s}{2s} (-\hat{x})^s \biggr)
  M_j
\end{equation}
eq.\ $(11)$ of \cite{Per}. Here
\begin{equation}\label{eq:3.7}
  \hat{x} M_j = M_{j-1}
\end{equation}
$\varepsilon_s$ for a graph $g$ is a linear function of a finite number of $g
\sslash \ell_i$, $\ell_i$ a set of given graphs, the 'contributors'.
The only thing we need to know is that for any given product of $\varepsilon_s$'s, $\prod_i
\varepsilon_{s(i)}$, one has that
\begin{equation}\label{eq:3.8}
  \EE\biggl( \prod_i \varepsilon_{s(i)} \biggr) \le C
\end{equation}
i.e.\ it is a bounded function of $n$. Here as everywhere in this paper the expectation
is the average value of the function
over all r regular bipartite graphs of order $2n$.

The result one needs to see this is that the number of $s$-cycles are
independent Poisson random variables of finite means in the fixed $r$, $n$ goes
to infinity limit, \cite{Bol}.
One then uses the fact that the $\ell_i$ and $g_{ij}$ graphs discussed above all
are either single cycle or multicycle in nature.

Working from eq.\ $(3.5)$ one can arrange the resultant terms arising 
into the following expression for $m_j$
\begin{align}
  m_j &= \frac{n^jr^j}{j!} (1+\hat{H_j)} \label{eq:3.2}\\
 \hat {H_j }&= \sum_{h=1}^{j-1} \frac{\hat{a_h}(r,j,\{\epsilon_i\})}{n^h} \label{eq:3.3}
\end{align}
$m_j$ is a function on graphs, eq.\ $(3.5)$ or eq.\ $(3.8)$-$(3.9)$ in turn expresses $m_j$ as a polynomial
in the $\{\epsilon_i\}$, these also functions on the graphs. We will be dealing with expectations of
polynomials in the $\{m_j\}$, for example eq.\ $(4.3)$. We make the important observation that,
for the sum in eq.\ $(3.9)$ appearing in an expectation,
the $n$ dependence of the $\epsilon_i$ does not effect the formal expected asymptotic expansion
by powers of $1/n$, from the discussion surrounding
eq.\ $(3.7)$. 

Assuming as we do the conjecture of Section $10$ there follows from eq.\ $(3.1)$-$(3.6)$ 
and  $(3.8)$-$(3.9)$ 

\begin{align}
   [j^k n^{-h}]\,\ln \biggl( 1 +\hat{ H_j }\biggr) =
    [j^k n^{-h}]\,\ln \biggl( 1 + \sum_{s=1}^{j-1} \frac{\hat{a_s}(r,j,\{\epsilon_i\})}{n^s} \biggr) &=
  0, \qquad k \ge h+2 \label{eq:3.4}\\
  [j^{h+1} n^{-h}]\,\ln \biggl( 1 +\hat{ H_j }\biggr)=  [j^{h+1} n^{-h}]\,\ln \biggl( 1 + \sum_{s=1}^{j-1} \frac{\hat{a_s}(r,j,\{\epsilon_i\})}{n^s} \biggr)
  &= \frac{1}{(h+1)h} \biggl( \frac{1}{r^h} - 2 \biggr) \label{eq:3.5}
\end{align}
Again we are working with formal polynomials, for fixed $r$ and $\epsilon_i$.

\section{Some simple reorganization}
We define
\begin{equation}\label{eq:4.1}
  1 + K_i \equiv \frac{(v-1)^i}{r^i} \cdot \frac{(v-2i)!\,i!\,2^i}{v!} \cdot
  \frac{r^in^i}{i!}
\end{equation}
using notably eq.\ \eqref{eq:1.2}.
Then with
\begin{equation}\label{eq:4.2}
  \alpha_0=\Biggl(
    \prod_{\ell \in \mathcal{L}^+} \bigl(
      (1+\hat{H_{i+\ell}})(1+K_{i+\ell})
    \bigr)^{\binom{k}{\ell}}
    - \prod_{\ell \in \mathcal{L}^-} \bigl(
      (1+\hat{H_{i+\ell}})(1+K_{i+\ell})
    \bigr)^{\binom{k}{\ell}}
  \Biggr)
\end{equation}
$\alpha$ becomes
\begin{equation}\label{eq:4.3}
  \alpha=\EE(\alpha_0)
 \end{equation}

Further we set
\begin{equation}\label{eq:4.3}
  1 + K_i \equiv e^{G_i}
\end{equation}
where
\begin{align}
  G_i &\equiv G_{i,1} + G_{i,2} + G_{i,3} + G_{i,4} + G_{i,5} \label{eq:4.4}\\
  G_{i,1} &\equiv i\,\ln\Bigl( 1 - \frac{1}{2n} \Bigr) \label{eq:4.5}\\
  G_{i,2} &\equiv (2n-2i)\,\ln\Bigl( 1 - \frac{i}{n} \Bigr) \label{eq:4.6}\\
  G_{i,3} &\equiv 2i \label{eq:4.7}\\
  G_{i,4} &\equiv \frac{1}{2}\,\ln\Bigl( 1 - \frac{i}{n} \Bigr) \label{eq:4.8}\\
  G_{i,5} &\equiv \sum_{j\ \text{odd}} c_j \biggl( \frac{1}{n^j} -
  \frac{1}{(n-i)^j} \biggr) \label{eq:4.9}
\end{align}
We have used the Stirling series to expand $\ln n!$.
We also note that for example $c_1 = -\frac{1}{24}$.
 $K_i$ is easily developed as a series in inverse powers of $n$.
 
 \textit{Valuable Observation} The convergence problem for series, except inside
 expectation values, is trivial, since one deals with $r$, $j$, and $\epsilon_i$ ( taken as a
 number ) fixed and $n$ large enough. BUT, the only expectations we take are of $\alpha_0$
 and $\alpha_0^{2}$ ( for $\beta$ ). And, see (4.2) and the discussion after (3.8)-(3.9), these
 both are finite polynomials in the $\{\epsilon_i\}$! ( So to study $\alpha_0$ and $\alpha$
 it is a good idea to expand $\alpha_0$ in the formal series
 in powers of $\frac{1}{n}$ taking the coefficients of the terms through $\frac{1}{n^{k-1}}$ from
 eq.(6.2) and the rest of the terms from eq.(4.2). )

\section{{$k = 1$} and {$k = 0$}} 
Not only is $k = 1$ the first case, but it is different from $k \ge 2$ in some
essential ways.
We proceed to compute $\alpha$ for $k = 1$.
From eq.\ \ref{eq:4.2} we have
\begin{equation}\label{eq:5.1}
  \alpha_0 = \bigl(
    (1+\hat{H_{i+1}})(1+K_{i+1}) - (1+\hat{H_{i})}(1+K_{i})
  \bigr)
\end{equation}
From eq.\ $(18)$ and eq.\ $(45)$ of \cite{Per} 
one gets
\begin{equation}\label{eq:5.2}
 \hat{ H_i} = i(i-1)\Bigl( -1 + \frac{1}{2r} \Bigr) \frac{1}{n} + \OO\Bigl(
  \frac{1}{n^2} \Bigr)
\end{equation}
In such asymptotic series bounds we treat the $\epsilon_i$ as constants.

Using eq.\ \eqref{eq:4.3}--\eqref{eq:4.9} one
gets
\begin{equation}\label{eq:5.3}
  K_i = (i^2 - i)\frac{1}{n} + \OO\Bigl( \frac{1}{n^2} \Bigr)
\end{equation}
There follows
\begin{thm}\label{thm:5.1}
  For $k = 1$
  \begin{equation}\label{eq:5.4}
    \alpha_0 = \frac{i}{rn} + \OO\Bigl( \frac{1}{n^2} \Bigr)
  \end{equation}
\end{thm}
One easily gets
\begin{thm}\label{thm:5.2}
  For $k = 0$
  \begin{equation}\label{eq:5.4a}
    \alpha_0= 1 + \OO\Bigl( \frac{1}{n} \Bigr)
  \end{equation}
\end{thm}

\section{$k \ge 2$}
The goal of this section is proving the following theorem.
\begin{thm}\label{thm:6.1}
  For $k \ge 2$
  \begin{equation}\label{eq:6.1}
 \alpha_0 = \frac{(k-2)!}{r^{k-1}n^{k-1}} + \OO(\frac{1}{n^{k}})
  \end{equation}
  
\end{thm}
From Section \ref{sec:8} using the Second Identity we have
\begin{equation}\label{eq:6.2}
  \alpha_0 = \biggl( \Bigl( 1 + t_+ + \frac{1}{2}t_+^2 + \cdots\Bigr) -
  \Bigl( 1 + t_- + \frac{1}{2}t_-^2 + \cdots\Bigr) \biggr)
\end{equation}
where with
\begin{equation}\label{eq:6.3}
  1 + U_i = (1 + \hat{H_i})(1 + K_i)\\
\end{equation}
one defines
\begin{equation}
 t_+ = \sum_{\ell \in \mathcal{L}^+} \binom{k}{\ell}
  \Bigl( U_{i+\ell} - \frac{1}{2} (U_{i+\ell})^2 + \frac{1}{3} (U_{i+\ell})^3
  \cdots \Bigr)\label{eq:6.4}\\
\end{equation}
\begin{equation}
  t_- = \sum_{\ell \in \mathcal{L}^-} \binom{k}{\ell}
  \Bigl( U_{i+\ell} - \frac{1}{2} (U_{i+\ell})^2 + \frac{1}{3} (U_{i+\ell})^3
  \cdots \Bigr)\label{eq:6.5}
\end{equation}
We treat the terms written explicitly in \eqref{eq:6.2}; the induction to
higher powers of $t$ is trivial.
\begin{align}
  \biggl[ \frac{1}{n^d} \biggr] (t_+ - t_-) &= \biggl[ \frac{1}{n^d} \biggr]
  \sum_\ell \binom{k}{\ell} (-1)^{k+\ell} \Bigl( U_{i+\ell} - \frac{1}{2}
  (U_{i+\ell})^2 + \frac{1}{3} (U_{i+\ell})^3 \cdots \Bigr)\label{eq:6.6}\\
  &= \begin{cases}
    \hfil 0 & d < k - 1\\
    \displaystyle\frac{(k-2)!}{r^{k-1}} & d = k-1
  \end{cases}\label{eq:6.7}
\end{align}
by the First Identity, Theorem \ref{thm:7.1}. In applying $\bigg[\frac{1}{n^{d}}\bigg]$
we treat the $\epsilon_i$ as constants.
Next we want to prove that the higher powers of $t$'s make no contribution in
\eqref{eq:6.2}!
This is amazing when one first sees it.
\par We want to show
\begin{equation}\label{eq:6.8}
  \biggl[ \frac{1}{n^d} \biggr] (t_+^2 - t_-^2) = 0 \quad \text{for} \quad d
  \le k - 1
\end{equation}
We proceed by looking at the powers of $\frac{1}{n}$.
\begin{equation}\label{eq:6.9}
  \biggl[ \frac{1}{n^d} \biggr] (t_+^2 - t_-^2) = \sum_{s=1}^{d-1} \Biggl[
    \biggl( \biggl[ \frac{1}{n^s} \biggr] t_+ \biggr)
    \biggl( \biggl[ \frac{1}{n^{d-s}} \biggr] t_+ \biggr)
  - \biggl( \biggl[ \frac{1}{n^s} \biggr] t_- \biggr)
    \biggl( \biggl[ \frac{1}{n^{d-s}} \biggr] t_- \biggr) \Biggr]
\end{equation}
All we need to complete a proof of \eqref{eq:6.8} is to show
\begin{equation}\label{eq:6.10}
  \biggl[ \frac{1}{n^s} \biggr] t_+ = \biggl[ \frac{1}{n^s} \biggr] t_- \qquad
  1 \le s \le d - 1
\end{equation}
But this follows from \eqref{eq:6.6},\eqref{eq:6.7} above.
Pretty neat.

\section{First Identity}
For convenience we introduce
\begin{equation}\label{eq:7.1}
  F_i \equiv (1 + \hat{H_i})(1 + K_i)
\end{equation}

\begin{thm}[First Identity]\label{thm:7.1}
  For all $r$, $i \ge 0$, and $k \ge 2$
  \begin{equation}\label{eq:7.2}
    \sum_{\ell=0}^k \binom{k}{\ell} (-1)^{\ell+k} \biggl[ \frac{1}{n^{k-1}}
    \biggr] \sum_{m=1}^{k-1} (-1)^{m+1} \frac{1}{m} \bigl( F_{(i+\ell)} - 1
    \bigr)^m = \frac{(k-2)!}{r^{k-1}}
  \end{equation}
\end{thm}
\begin{thm}\label{thm:7.2}
  For all $r$, and $k \ge 2$
  \begin{equation}\label{eq:7.3}
    \biggl[ \frac{1}{n^{k-1}} \biggr] \ln(F_i)
  \end{equation}
  has highest power of $i = i^k$, and this term is
  \begin{equation}\label{eq:7.4}
    \frac{(k-2)!}{k!}\frac{i^k}{r^{k-1}}
  \end{equation}
\end{thm}
For example for $k = 3$
\begin{equation}\label{eq:7.5}
  \biggl[\frac{1}{n^2}\biggr] \ln(F_s) = -\frac{1}{12} \frac{s\bigl(3r^2s - 3r^2
  - 12rs - 2s^2 + 12r + 9s - 7\bigr)}{r^2}
\end{equation}
\par We now note
\begin{equation}\label{eq:7.6}
  \sum_{\ell=0}^k \binom{k}{\ell} (-1)^{\ell+k} \ell^d = \begin{cases}
    \hfil 0 & d < k\\
    k! & d = k
  \end{cases}
\end{equation}
that follows from
\begin{equation}\label{eq:7.7}
  \sum_{\ell=0}^k \binom{k}{\ell} (-1)^{\ell+k} \ell^d = \Delta^k i^d
\end{equation}
which like $(\frac{d}{dx})^k x^d$ has values in \eqref{eq:7.6}.
From eq.\ \eqref{eq:7.6} one can deduce that Theorem \ref{thm:7.1} follows from
Theorem \ref{thm:7.2}, which we proceed to prove.
\begin{equation}\label{eq:7.8}
  \biggl[ \frac{1}{n^{k-1}} \biggr] \ln(F_i) = 
  \biggl[ \frac{1}{n^{k-1}} \biggr] \ln(1+\hat{H_i})
  + \biggl[ \frac{1}{n^{k-1}} \biggr] \ln(1+K_i)
\end{equation}
From eqs.(3.10) (3.11) with $k \ge 2$ we see that the highest
power of $i$ in $[\frac{1}{n^{k-1}}]\ln(1+\hat{H_i)}$ is $k$ and its coefficient is
\begin{equation}\label{eq:7.9}
  \frac{1}{k(k-1)} \biggl( \frac{1}{r^{k-1}} - 2 \biggr)
\end{equation}
To study $[\frac{1}{n^{k-1}}]\ln(1+K_i)$ we turn to equations
(4.4)-(4.10).
We note the highest power of $i$ arises from the expansion of the term
$G_{i,2}$, eq.(4.7), and $[\frac{1}{n^{k-1}}]\ln(1+K_i)$ has highest
power $i$
\begin{equation*}
  \frac{2}{k(k-1)}
\end{equation*}
So
\begin{equation}\label{eq:7.10}
  \biggl[ \frac{1}{n^{k-1}} \biggr] \ln(F_i) = \frac{(k-2)!}{k!}
  \frac{1}{r^{k-1}} i^k
\end{equation}
Quod erat demonstrandum.

\section{Second Identity}\label{sec:8}
We start with some simple manipulations
\begin{equation}\label{eq:8.1}
  \prod_i (1+x_i)^{e_i} = e^{\sum e_i\ln(1+x_i)} = 1 + \Bigl( \sum e_i\ln(1+x_i)
  \Bigr) + \frac{1}{2!} \Bigl( \sum e_i\ln(1+x_i) \Bigr)^2 + \cdots
\end{equation}
With the notation
\begin{equation}\label{eq:8.2}
  (1 + \hat{H_i)}(1 + K_i) \equiv 1 + U_i
\end{equation}
we substitute $U_i$ for $x_i$ and $\binom{k}{\ell}$ for $e_i$ in \eqref{eq:8.1}
\begin{equation}\label{eq:8.3}
  \prod_{\ell \in \mathcal{L}^+} (1 + U_{i+\ell})^{\binom{k}{\ell}} = 1 + t _+ +
  \frac{1}{2} t_+^2 + \frac{1}{3!} t_+^3 \cdots
\end{equation}
where
\begin{equation}\label{eq:8.4}
  t _+\equiv \sum_{\ell \in \mathcal{L}^+} \binom{k}{\ell} \biggl( U_{i+\ell} -
  \frac{1}{2} (U_{i+\ell})^2 + \frac{1}{3} (U_{i+\ell})^3 \cdots \biggr)
\end{equation}
The Second Identity consists of \eqref{eq:8.3} and \eqref{eq:8.4} and the same
expressions with $\mathcal{L}^+$, $t_+$ replaced by $\mathcal{L}^-$, $t_-$.

\section{Completion}

The information we need from the calculations of this paper are \textbf{ Theorem 5.1},
\textbf{Theorem 5.2}, and \textbf{Theorem 6.1}.  From these respectively we get:

1) For $k=1, i \neq 0$
\begin{align}\label{eq:9.1}
  \alpha &\ge \frac{c}{n} ,  \qquad    c \quad  positive\\
 \beta &\le \frac{c}{n^{4}} ,  
\end{align}

2) For $k=0$
\begin{align}\label{eq:9.2}
  \alpha &\ge c ,  \qquad    c \quad  positive\\
 \beta &\le \frac{c}{n^{2}} ,  
\end{align}

3) For $k \ge 2$
\begin{align}\label{eq:9.3}
  \alpha &\ge \frac{c}{n^{k-1}} ,  \qquad    c \quad  positive\\
 \beta &\le \frac{c}{n^{2k}} ,  
\end{align}

We are notationally using a $c$ that varies from equation to equation. 
( We due not pursue stronger results that follow from the fact that $\epsilon_i$
is zero for $i=1,2,3$ among other possible improvements. ) Referring to Section 2,
\textbf{Theorem 1.1} follows from fact that $\frac{\beta}{\alpha^2}$ goes to zero as $n$ goes 
to infinity in each case.

\section{An Awesome Conjecture}
Let $z_i$ be positive integers. We set:
\begin{equation}\label{eq:10.1}
F =  \sum_{s \ge 0} \frac{a_s(r,j)}{n^s} +\sum _i c_i j (j-1)\cdots (j-z_i+1)\frac{1}{n^{z_i} r^{z_i}} \sum_{s \ge 0} \frac{a_s(r,j-z_i)}{n^s} 
\end{equation}
Then we conjecture:
\begin{align}
  [j^k n^{-h}]\,\ln ( F ) &=
  0, \qquad k \ge h+2 \label{eq:10.2}\\
  [j^{h+1} n^{-h}]\,\ln ( F )
  &= \frac{1}{(h+1)h} \biggl( \frac{1}{r^h} - 2 \biggr) \label{eq:10.3}
\end{align}
Compare eq. (10.2) - (10.3) to eq. (3.4) - (3.5).

\textbf{Acknowledgement}   We thank Ian Wanless for enlightening us on a number of points.

\end{document}